\documentclass[12pt]{amsart}

\usepackage{amsfonts,amsmath,amssymb,amsthm}
\usepackage{color}
\usepackage{marginnote}

\usepackage{lipsum}

\newtheorem{theorem}{Theorem}
\theoremstyle{plain}

\newtheorem{corollary}{Corollary}

\newtheorem{definition}{Definition}

\newtheorem{proposition}{Proposition}
\newtheorem{remark}{Remark}

\numberwithin{equation}{section}

\DeclareMathOperator{\Ric}{Ric}
\DeclareMathOperator{\grad}{grad}

\DeclareMathOperator{\diver}{div}

\begin{document}
\thispagestyle{empty}
\title[Null Sectional Curvatures of Warped Products]
{Null Sectional Curvatures of Warped Products}
\date{July 26, 2016}
\subjclass[2010]{Primary 53C25; Secondary 53C80}
\keywords{Warped products, null sectional curvature,
warped product space-times.}

\begin{abstract}
In this paper, we investigate the null (light-like) sectional curvatures
of Lorentzian warped product manifolds. We derive the formulas
for the null sectional curvature of many well-known
warped product space-time models such as multiply generalized
Robertson-Walker space-times, generalized Kasner space-times and
standard static space-times.
\end{abstract}

\author{Beng\.{I} R. Yavuz}
\address[B. R. Yavuz]{Department of Mathematics, Bilkent University,
         Bilkent, 06800 Ankara, Turkey}
\email{bengi@fen.bilkent.edu.tr}
\author{B\"{u}lent \"{U}nal}
\address[B. \"{U}nal]{Department of Mathematics, \!Bilkent University,
         \!Bilkent, 06800 \!\! Ankara, Turkey}
\email{bulentunal@mail.com}
\author{Fernando Dobarro}
 \address[F. Dobarro]{
 Instituto de Desarrollo Econ\'{o}mico e Innovaci\'{o}n, Universidad
 Nacional de Tierra del Fuego, Ant\'{a}rtida e Islas del
 Atl\'{a}nti-co Sur, Gobernador Paz 1410, 9410, Ushuaia, Tierra del
 Fuego, Argentina.
 \newline
 \hbox{\hspace{2.3 cm}} FCE, Universidad Austral, Mariano Acosta s/n
 y Ruta Nac. 8, Edificio Grado, B1629WWA, Buenos Aires, Argentina}
 \email{fdobarro@untdf.edu.ar \\ fdob07@gmail.com}

\pagenumbering{gobble}
\clearpage \thispagestyle{empty}

\maketitle

\tableofcontents

\newpage

\renewcommand{\thepage}{\arabic{page}}
\setcounter{page}{1}

\section{Introduction}

Warped product manifolds were first introduced to the literature
by R. Bishop and B. O'Neill in \cite{BishopON} to construct
complete Riemannian manifolds with negative sectional curvature
everywhere. After that Beem, Ehrlich and Powell established that
a wide class of well known exact solutions of the Einstein's
field equations can be expressed as a Lorentzian warped product in
\cite{BEP1982, BEE1996}. Two important examples of warped products
are the generalized Robertson-Walker
space-times (GRW, for short) and the standard static space-times
(SSS-T, for short). The former are obviously a generalization of
Robertson-Walker space-times and the latter a generalization of the
Einstein static universe. In addition to these space-times, we can list
generalized Kasner and Reissner-Nordstr\"om space-time models.

In the current paper, we focus on the null (light-like) sectional
curvatures of Lorenztian warped products. The concept of null sectional
curvature were first defined by S. G. Harris in \cite{Harris} to
study sectional curvatures of null plane sections of Lorentzian manifolds.
In \cite{Harris85}, it is proved that Robertson-Walker metrics can be locally
characterized as those for which the null sectional curvature denoted by $K_{N}$
at each point is a constant for all the null planes at that point.

In \cite{Palomo2007}, J. P. Palomo proved the smoothness
of the null sectional curvature. Moreover, in a not necessarily complete
manifold, if the null sectional curvature is not zero for any degenerate
planes at a point $p$, then there is at most one local GRW structure in a
neighborhood of $p$ (see \cite{GO}). Moreover he further showed that for
an $n(\geq 4)$-dimensional manifold the $U-$ normalized null sectional curvature
function is constant on every pencil of degenerate planes
belonging to a certain null direction, if and only if, the
Lorentzian manifold is conformally flat.

In \cite{Harris85, Koch}, it is shown that if
$(M,g)$ is a time orientable Lorentzian manifold with
$dim(M) \geq 3,$ and $U$ is a globally defined unitary time-like
vector field, then the $U-$ normalized null sectional curvature
is isotropic, i.e, it is only a point function.
Moreover,

$$K^U(p, \Pi) = K^U(p)$$

for all null planes $\Pi$ if and only if the curvature tensor
satisfies:

(a) $R(X,Y)Z = k (X \wedge_g Y)Z,$ for any $X,Y,Z \in U^{\perp}$

(b) $R(X,U)U = \mu X,$ for any with
$\kappa , \mu \colon M \to \mathbb R$ and $K^U = \kappa + \mu.$

If the $U-$ normalized null sectional curvature $K^U$ is isotropic
(it is said to be spatially constant) if it is constant on the space
of orthogonal to the chosen time-like vector field $U$, i.e, $K^U=0$
for every $X \in U^{\perp}$ (see \cite{Koch}).

Let $(M,g)$ be an $n \geq 4$- dimensional Lorentzian manifold
and let $U$ be a unitary time-like vector field on $M.$ Suppose
that the $U$-normalized null sectional curvature is non-zero,
isotropic and spatially constant. Then $g$ is locally a
Friedmann-Lemaitre-Robertson-Walker metric (see \cite{Harris85, Koch}).

In \cite{GO}, to characterize global decomposition of a manifold
as a generalized Robertson-Walker space-time, the null sectional
curvature is applied.

The null sectional curvature has been used in the study of
conjugate points along null geodesics (see \cite{GO1, Harris, Palomo2006}).

We recall that a \textit{warped product} can be defined as follows
\cite{BEE1996,ON}. Let $(B,g_B)$ and $(F,g_F)$ be
pseudo-Riemannian manifolds and also let $b \colon B \to
(0,+\infty)$ be a smooth function. Then the (singly) warped product,
$B \times {}_b F$ is the product manifold
$B \times F$ furnished with the metric tensor $g=g_B \oplus
b^{2}g_F$, more precisely
\begin{equation*}
g=\pi^{\ast}(g_B) \oplus (b \circ \pi)^2 \sigma^{\ast}(g_F),
\end{equation*}
 where $\pi \colon B \times F
\to B$ and $\sigma \colon B \times F \to F$ are the usual
projection maps and ${}^\ast$~denotes the pull-back
operator on tensors.

A \textit{standard static space-time} can be
considered as a Lorentzian warped product where the warping
function is defined on a Riemannian manifold called the
fiber and acting on the negative definite metric on an open
interval of real numbers, called the base. More precisely,
a SSS-T, ${}_f(a,b) \times F$ is a
Lorentzian warped product furnished with the metric
$g=-f^2{\rm d}t^2 \oplus g_F,$ where $(F,g_F)$ is a
Riemannian manifold, $f \colon F \to (0,+\infty)$ is smooth,
and $-\infty \leq a < b \leq +\infty.$ In \cite{ON}, it was
shown that any static space-time is locally isometric to a
SSS-T.

Standard static space-times have been previously studied by
many authors. Kobayashi and Obata \cite{KO} stated the
geodesic equation for this class of space-times and the
causal structure and geodesic completeness was considered
in \cite{AD}, where sufficient conditions on the warping
function for nonspacelike geodesic completeness of the
SSS-T was obtained (see also
\cite{RASM}). In \cite{AD1}, conditions are found which
guarantee that SSS-Ts either satisfy
or else fail to satisfy certain curvature conditions from
general relativity. The existence of geodesics in SSS-Ts
have been studied by several authors.
S\'{a}nchez \cite{mS2} gives a good overview of geodesic
connectedness in semi-Riemannian manifolds, including a
discussion for SSS-Ts.

Two of the most famous examples of SSS-Ts are Minkowski
space-times and the Einstein
static universe \cite{BEE1996,HE} which is $\mathbb R \times
\mathbb S^3$ equipped with the metric
$$g=-{\rm d}t^2+({\rm d}r^2+\sin^2 r {\rm d}\theta^2+
\sin^2 r \sin^2 \theta {\rm d} \phi^2)$$ where $\mathbb
S^3$ is the usual 3-dimensional Euclidean sphere and the
warping function $f \equiv 1.$ Another well-known example
is the universal covering space of the anti-de Sitter
space-time, a SSS-T of the form ${}_f
\mathbb R \times \mathbb H^3$ where $\mathbb H^3$ is the
3-dimensional hyperbolic space with constant negative
sectional curvature and the warping function $f \colon
\mathbb H^3 \to (0,+\infty)$ defined as
$f(r,\theta,\phi)=\cosh r$  \cite{BEE1996,HE}. Finally, we can
also mention the Exterior Schwarzschild space-time
\cite{BEE1996,HE}, a SSS-T of the form
${}_f \mathbb R \times (2m,+\infty) \times \mathbb S^2,$
where $\mathbb S^2$ is the 2-dimensional Euclidean sphere,
the warping function $f \colon (2m,+\infty) \times \mathbb
S^2 \to (0,+\infty)$ is given by
$f(r,\theta,\phi)=\sqrt{1-2m/r}$ and the line
element on $(2m,+\infty) \times \mathbb S^2$ is $${\rm
d}s^2=\bigl(1-\frac{2m}{r} \bigl)^{-1} {\rm d}r^2+r^2({\rm
d} \theta^2+\sin^2 \theta {\rm d} \phi^2).$$

\section{Preliminaries}

In this section, we give the formal definitions of several
types of warped product space-time models such as
GRWs, SSS-Ts
and multiply generalized Robertson-Walker space-times (MGRW, for short)
and their related geometric structure formulas (see \cite{BEE1996,ON}).

Throughout the article $I$ will denote an open real
interval of the form $I=(t_1,t_2)$, where $-\infty \leq t_1 < t_2
\leq +\infty$.

\begin{definition} \label{dfnr} Let $(F,g_F)$ be an $s$-dimensional
Riemannian manifold and $b \colon I \to (0,+\infty)$ be a
smooth function. Then the $n(=1+s)$-dimensional product
manifold $I \times F$ furnished with the metric tensor
$g=-{\rm d}t^2 \oplus b^2 g_F$ is called a \textit{generalized
Robertson-Walker space-time} and is denoted by $ I \times _b F$, where
${\rm d}t^2$ is the Euclidean metric tensor on $I.$
\end{definition}

\begin{definition} \label{dfna} Let $(F,g_F)$ be an $s$-dimensional
Riemannian manifold and $f \colon F \to (0,+\infty)$ be a
smooth function. Then the $n(=1+s)$-dimensional product
manifold $I \times F$ furnished with the metric tensor
$g=-f^2{\rm d}t^2 \oplus g_F$, where
${\rm d}t^2$ is the Euclidean metric tensor on $I$, is called a \textit{standard static
space-time} and is denoted by $_f I \times F$.
\end{definition}

\begin{definition}Let $(B,g_B)$ and $(F_i,g_{F_i})$ be pseudo-Riemannian
manifolds and also let $b_i \colon B \to (0,\infty)$ be smooth
functions for any $i \in \{1,2,\cdots,m\}.$ The product manifold $M=B \times F_1
\times F_2 \times \cdots \times F_m$ furnished with the metric
tensor $g=g_B \oplus b_{1}^{2}g_{F_1} \oplus b_{2}^{2}g_{F_2}
\oplus \cdots \oplus b_{m}^{2}g_{F_m}$ is called a {\it multiply
warped product} and is denoted as
\begin{equation*}
B \times {}_{b_1} F_1 \times \dots \times {}_{b_m} F_m.
\end{equation*}
More precisely
\begin{equation} \label{idwp}
g=\pi^{\ast}(g_B) \oplus (b_1 \circ \pi)^2
\sigma_{1}^{\ast}(g_{F_1}) \oplus \cdots \oplus (b_m \circ \pi)^2
\sigma_{m}^{\ast}(g_{F_m}),
\end{equation}
where the $\pi$ and $\sigma_{i}$ are the usual
projection maps and ${}^\ast$~denotes the pull-back
operator on tensors.
$(B,g_B)$, $b_i \colon B \to (0,\infty)$ and $(F_i,g_{F_i})$ are
called \textit{base}, $i$-\textit{fiber} and
$i$-\textit{warping function} of the multiply warped product, respectively.
\end{definition}

\begin{itemize}
\item If $m=1,$ then we obtain a {\it singly warped product}.

\item If all $b_i \equiv 1,$ then we have a (trivial) {\it product
manifold}.

\item If $(B,g_B)$ and all the $(F_i,g_{F_i})$ are {\it Riemannian}
manifolds, then $(M,g)$ is also a {\it Riemannian} manifold.

\item The multiply
warped product $(M,g)$ is a {\it Lorentzian multiply warped
product} if $(F_i,g_{F_i})$ are all {\it Riemannian} for any $i
\in \{1,2,\cdots,m\}$ and either $(B,g_B)$ is {\it Lorentzian} or
else $(B,g_B)$ is a one-dimensional manifold with a {\it negative
definite} metric $-{\rm d}t^2$.

\item If $B$ is an open
interval $I$ equipped with the negative
definite metric $g_B=-{\rm d}t^2$ and all the $(F_i,g_{F_i})$ manifolds
are Riemannian, then the Lorentzian multiply warped product
$(M,g)$ is called a \textit{multiply generalized Robertson-Walker
space-time}. In particular, a MGRW is called a \textit{generalized
Reissner-Nordstr\"om space-time} when $m=2.$
\end{itemize}

Throughout the paper the fiber(s) $(F,g_F)$ of a warped product
space-time model is always assumed to be connected. We denote the
set of lifts of vector fields on $B$ and $F$ to $B \times F$ by
$\mathfrak L(B)$
and $\mathfrak L(F),$ respectively and use the same notation for
a vector field and for its lift (see page
205 of \cite{ON}).

\noindent From now on, we follow the convention applied in
\cite{BEE1996} (note the difference with \cite{ON}) for
the definition and sign of the Riemann curvature tensor $R$,
namely. For any $n$-dimensional pseudo-Riemannian manifold $(N,h)$,
$$
R(X,Y)Z=\nabla_X \nabla_Y Z - \nabla_Y \nabla_X Z - \nabla_{[X,Y]} Z,
$$
where $\nabla$ is the $h$-Levi-Civita connection and $X, Y, Z$
are vector fields on $N$.
\noindent Besides, ${\rm Ric}$ will denote the Ricci tensor (see
\cite{ON}).

\noindent Furthermore, we will apply the sign
convention for the Laplacian in \cite{ON}, i.e,
$\Delta={\rm tr}({\rm H})=\diver \grad $, where
${\rm tr}$ denotes the $h$-trace, $H$ the Hessian tensor respect
to the Levi-Civita connection, $\diver $ the $h$-divergence and
$\grad $ the $h$-gradient (see page 85 of \cite{ON}).

\noindent Eventually, we will put a sub- or super-index in
the corresponding operator indicating the manifold on which
it is acting, for instance $\nabla ^{B}$ if the connection is
that on the manifold $B$.

The following basic formulas about the geometry of multiply warped
product can be found in \cite{Dobarro-Unal2005}

\begin{proposition} \label{gcovd} Let
$M=B \times {}_{b_1}F_1 \times \cdots \times {}_{b_m}F_m$ be a
{\it pseudo-Riemannian} multiply warped product with metric $g=g_B
\oplus b_{1}^{2}g_{F_1} \oplus \cdots \oplus b_{m}^{2}g_{F_m}$
also let  $X,Y \in \mathfrak L(B)$ and $V \in \mathfrak L(F_i),$
$W \in \mathfrak L(F_j).$ Then
\begin{enumerate}
\item ${\displaystyle \nabla_{X} Y= \nabla_{X}^B Y}$ \item
${\displaystyle \nabla_{X} V= \nabla_{V} X=\frac{X(b_i)}{b_i} V }$
\item $ \nabla_{V} W=
  \begin{cases}
{\displaystyle 0}&
\text{if $ i \neq j $ }, \\
{\displaystyle \nabla_{V}^{F_i} W-\frac{g(V,W)}{b_i}
\grad_{B}b_i } & \text{if $ i=j, $ }
\end{cases} \nonumber $
\end{enumerate}
where $\nabla = \nabla^M$.
\end{proposition}

One can compute the {\it gradient} and the {\it Laplace-Beltrami}
operator on $M$ in terms of the {\it gradient} and the {\it
Laplace-Beltrami} operator on $B$ and $F_i,$ respectively.

\begin{proposition} \label{ggra-lap}
Let $M=B \times {}_{b_1}F_1 \times \cdots \times {}_{b_m}F_m$ be a
{\it pseudo-Riemannian} multiply warped product with metric $g=g_B
\oplus b_{1}^{2}g_{F_1} \oplus \cdots \oplus b_{m}^{2}g_{F_m}$ and
$\phi \colon B \to \mathbb R$ and ${\psi}_i \colon F_i \to \mathbb
R$ be smooth functions for any $i \in \{1,\cdots,m\}.$ Then
\begin{enumerate}
\item $\grad(\phi \circ \pi)=\displaystyle{\grad_{B}\phi}$

\item $\grad({\psi}_i \circ {\sigma}_i)=\displaystyle{
\frac{1}{b_{i}^2}}\grad_{F_i}{\psi}_i$

\item $\Delta(\phi \circ \pi)=\displaystyle{\Delta_{B}\phi+
\sum_{i=1}^m s_i \frac{g_{B}(\grad_{B}\phi,\grad_{B}b_i)}{b_i}}$

\item $\Delta({\psi}_i \circ {\sigma}_i) =\displaystyle{
\frac{\Delta_{F_i}{\psi}_i}{b_{i}^2}}$
\end{enumerate}
where $\Delta = \Delta_M$ and $\grad = \grad_M$.
\end{proposition}

\begin{proposition}\label{grcur} Let $M=B \times {}_{b_1}F_1 \times \cdots
\times {}_{b_m}F_m$ be a {\it pseudo-Riemannian} multiply warped
product with metric $g=g_B \oplus b_{1}^{2}g_{F_1} \oplus \cdots
\oplus b_{m}^{2}g_{F_m}$ also let $X,Y,Z \in \mathfrak L(B)$ and
$V  \in \mathfrak L(F_i), W  \in \mathfrak L(F_j)$ and $U \in
\mathfrak L(F_k).$ Then
\begin{enumerate}

\item $\displaystyle{R(X,Y)Z = R_{B}(X,Y)Z}$

\item $\displaystyle{R(V,X)Y = -\frac{{\rm H}_B^{b_i}(X,Y)}{b_i}V }$

\item $\displaystyle{R(X,V)W=R(V,W)X=R(V,X)W=0}$ if \, $i \neq j.$

\item $\displaystyle{R(X,Y)V=0}$

\item $\displaystyle{R(V,W)X=0}$ if \, $i=j.$

\item $\displaystyle{R(V,W)U=0}$ if \, $i=j$ and $i,j
\neq k.$

\item $\displaystyle{R(U,V)W=-g(V,W)\frac{g_{B}(\grad_{B}b_i,
\grad_{B}b_k)}{b_i b_k}U}$ if \, $i=j$ and $i,j \neq k.$

\item $\displaystyle{R(X,V)W=-\frac{g(V,W)}{b_i}\nabla^{B}_{X}
(\grad_{B}b_i)}$ if \, $i=j.$

\item $R(V,W)U=
R_{F_i}(V,W)U+ \\
\displaystyle{\frac{\| \grad_{B}b_i \|_{B}^{2}} {b_{i}^2}
\left(g(V,U)W-g(W,U)V \right)} \textit{ if } i,j=k.$
\end{enumerate}
where $R = R_M$.
\end{proposition}

\begin{proposition}\label{gricur} Let $M=B \times {}_{b_1}F_1 \times \cdots
\times {}_{b_m}F_m$ be a {\it pseudo-Riemannian} multiply warped
product with metric $g=g_B \oplus b_{1}^{2}g_{F_1} \oplus \cdots
\oplus b_{m}^{2}g_{F_m},$ also let $X,Y,Z \in \mathfrak L(B)$ and
$V  \in \mathfrak L(F_i)$ and $W  \in \mathfrak L(F_j).$ Then
\begin{enumerate}
\item $\displaystyle{\Ric(X,Y)=\Ric_{B}(X,Y)-\sum_{i=1}^m
\frac{s_i}{b_i}{\rm H}_{B}^{b_i}(X,Y)}$ \item
$\displaystyle{\Ric(X,V)=0}$
\item $\displaystyle{\Ric(V,W)=0}$ if
\, $i \neq j.$
\item $\displaystyle{\Ric(V,W)=\Ric_{F_i}(V,W)-
\Big(\frac{\Delta_{B}b_i}{b_i}+(s_i-1) \frac{\| \grad_{B}b_i
\|_{B}^{2}}{b_{i}^2}}  \\
\displaystyle{+\sum_{
{k=1,
k \neq i }}^m s_k \frac{g_{B}(\grad_{B}b_i,\grad_{B}b_k)}{b_i
b_k}\Bigl)g(V,W)}$ if \, $i=j,$
\end{enumerate}
where $\Ric = \Ric_M$.
\end{proposition}

\section{Null Sectional Curvature}

Suppose that $\Pi$ is a null plane, that is $\Pi$ consists of a
one-dimensional subspace of null vectors and of space-like
vectors perpendicular to that subspace. Let $L$ be one
of the null vectors and $S$ be one of the space-like vectors.
Since $Q(\Pi)=0$, sectional curvature is not defined for a null plane.
Then S. Harris introduces the {\it null sectional curvature}
in \cite{Harris} for degenerate planes as follows
(see also \cite{AlbujerHaesen2010}):

In order to define {\it null sectional curvature},
first we need to fix a choice of null vector.

\begin{definition}
Let $M$ be an $n$-dimensional Lorentzian manifold.
A {\it null congruence} on $M$ is a submanifold $C$ of
the tangent bundle ${\rm T}_{0}M$ of nonzero null vectors on
$M$ such that for all $N$ in ${\rm T}_{0}M$,
there is exactly one scalar $\alpha$ satisfying $\alpha N \in C$.
\end{definition}

Then, {\it null sectional curvature} of $\Pi$ with respect to $N$ is
given by
$${ \rm K}_{N}(\Pi)=\frac{g(R(L,S)S,L)}{g(S,S)}  $$
where $R$ is the Riemannian curvature tensor.

\begin{remark} Note that,
this formula is independent of the choice of the space-like
vector $S$ in $\Pi$. However,
$K_N(p, \Pi) $
depends quadratically on the null vector $L$.
\end{remark}

Null sectional curvature can be normalized via the help of a
time-like vector field $U$ in the following way:
Assume that $U$ is a time-like vector field on a Lorentzian
manifold $M.$ Then the null congruence
$C(U)$ associated with $U$ is given by
$$C(U)=\{L \in \rm T_{0}M | g(L,L)=0, g(L,U)=-1 \}.$$

\begin{remark}The null congruence lies in the future null cone
due to the -1 in the definition.
\end{remark}

Then $U$-normalized null sectional curvature is

$${ \rm K}_{U}(\Pi)=\frac{g(R(L,S)S,L)}{g(S,S)}$$
where
$$C_U(M) =\{L \in TM | \, g(L,L)=0 \quad \text{and}
\quad g(L, U_{\Pi(L)}) = -1 \}.$$

\subsection{Null Sectional Curvature of a multiply GRW space-time}

Let $M=I \times_{b_{1}}  F_{1} \ldots \times_{b_{m}}  F_{m}$
be a multiply GRW space-time with the metric
$g=-{ \rm d}t^{2} \oplus b_{1}^{2}g_{F_{1}} \oplus \ldots \oplus b_{m}^{2} g_{F_{m}}$.
Let $L=- \partial_{t} +\mathbf{V} $ and $S=Y+\mathbf{W}$ where
$\mathbf{V} =\mathbf{\Sigma} V_{i}$    $\mathbf{W} = \mathbf{\Sigma} W_{j}$.
Then,

\begin{equation*} 
\begin{split}
g(R(L,S)S,L) &= g(R(-\partial_{t} + \mathbf{V},Y+\mathbf{W})\mathbf{V},
Y+\mathbf{W},-\partial_{t} + \mathbf{V}) \\
&= g(R(-\partial_{t},Y)Y,-\partial_{t}) +
g(R(\mathbf{\Sigma} V_{i},Y)Y,-\partial_{t}) \\
&+ g(R(-\partial_{t},\mathbf{\Sigma} W_{j})Y,-\partial_{t}) +
g(R(\mathbf{\Sigma} V_{i},\mathbf{\Sigma} W_{j})Y,-\partial_{t}) \\
&+ g(R(-\partial_{t},Y)\mathbf{\Sigma} W_{j},-\partial_{t})+
g(R(\mathbf{\Sigma} V_{i},Y)\mathbf{\Sigma} W_{j},-\partial_{t}) \\
&+  g(R(-\partial_{t},\mathbf{\Sigma} W_{j})
\mathbf{\Sigma} W_{j},-\partial_{t}) + g(R(\mathbf{\Sigma} V_{i},
\mathbf{\Sigma} W_{j})\mathbf{\Sigma} W_{j},-\partial_{t}) \\
&+ g(R(-\partial_{t},Y)Y,\mathbf{\Sigma} V_{i}) + g(R(\mathbf{\Sigma}
V_{i},Y)Y,\mathbf{\Sigma} V_{i}) \\
&+ g(R(-\partial_{t},\mathbf{\Sigma} W_{j})Y,\mathbf{\Sigma}
V_{i}) + g(R(\mathbf{\Sigma} V_{i},\mathbf{\Sigma} W_{j})Y,
\mathbf{\Sigma} V_{i}) \\
&+ g(R(-\partial_{t},Y)\mathbf{\Sigma} W_{j},\mathbf{\Sigma}
V_{i}) + g(R(\mathbf{\Sigma} V_{i},Y)\mathbf{\Sigma} W_{j},
\mathbf{\Sigma} V_{i}) \\
&+ g(R(-\partial_{t},\mathbf{\Sigma} W_{j})\mathbf{\Sigma}
W_{j},\mathbf{\Sigma} V_{i})+ g(R(\mathbf{\Sigma} V_{i},
\mathbf{\Sigma} W_{j})\mathbf{\Sigma} W_{j},\mathbf{\Sigma} V_{i})
\end{split}
\end{equation*}

\begin{equation} \label{mgrwsect}
\begin{split}
g(R( \mathbf{V},\mathbf{W})\mathbf{W}),\mathbf{V}) &=
\mathbf{\Sigma}_{i \neq j} -\frac
{1}{b_{j}^2}b_{j}^{\prime \prime}g(W_{j},W_{j})  g( V_{i},V_{i}) \\
&+ \mathbf{\Sigma} g(R_{F_{i}}(V_{i},W_{i}),W_{i},V_{i}) \\
&+ \mathbf{\Sigma}\frac{1}{b_{i}^2}b_{i}^{\prime \prime} g(g(V_{i},W_{i})W_{i}-g(W_{i},W_{i})V_{i}),V_{i}) \\
&= \mathbf{\Sigma}_{i \neq j} -(b_{j}^{\prime})^2
g_{F_{j}}(W_{j},W_{j})b_{i}^2g_{F_{i}}(V_{i},V_{i}) \\
&+ \mathbf{\Sigma} b_{i}^2g_{F_{i}}(R_{F_{i}}(V_{i},W_{i}),W_{i},V_{i}) \\
&+ \mathbf{\Sigma} b_{i}^2 (b_{i}^{\prime})^2[g_{F_{i}}(V_{i},W_{i})^2-
g_{F_{i}}(V_{i},V_{i})g_{F_{i}}(W_{i},W_{i})]
\end{split}
\end{equation}

\begin{theorem} \label{athm}
Let $M=I \times_{b_{1}}  F_{1} \ldots \times_{b_{m}}  F_{m}$
be a multiply GRW space-time with the metric
$g=-{ \rm d}t^{2} \oplus b_{1}^{2}g_{F_{1}} \oplus \ldots \oplus b_{m}^{2} g_{F_{m}}$.
Assume that $L=- \partial_{t} +\mathbf{V} $ and $S=Y+\mathbf{W}$ where
$\mathbf{V} =\mathbf{\Sigma} V_{i}$ and $\mathbf{W} = \mathbf{\Sigma} W_{j}$.
Then null sectional of $(M,g)$ is given by
\begin{equation*}
\begin{split}
{ \rm P}(p,\Pi) &=  \mathbf{\Sigma} b_{k}g_{F_{k}}(W_{k},V_{k})H^{B}_{b_{k}}(\partial_{t},Y) \\
&+ \mathbf{\Sigma} b_{k}g_{F_{k}}(V_{k},V_{k}){ \rm H}^{B}_{b_{k}}(Y,Y) \\
&- \mathbf{\Sigma} b_{j} b_{j}^{\prime \prime}g_{F_{j}}(W_{j},W_{j}) \\
&+ \mathbf{\Sigma} b_{i}g_{F_{i}}(W_{i},V_{i}){ \rm H}_{B}^{b_{i}}(\partial_{t},Y) \\
&- \mathbf{\Sigma}_{j \neq k} b_{k} (b_{j}^{\prime})^{2}g_{F_{k}}(V_{k},V_{k}) g_{F_{j}}(W_{j},W_{j}) \\
&+ \mathbf{\Sigma} b_{i}^{2} g_{F_{i}}(R_{F_{i}}(W_{i},V_{i})V_{i},W_{i}) \\
&- \mathbf{\Sigma} b_{i}^{2}(b_{i}^{\prime})^{2} b_{i}^{\prime \prime}[(g_{F_{i}}(W_{i},V_{i}))^{2}
-g_{F_{i}}(V_{i},V_{i})g_{F_{i}}(W_{i},W_{i})]
\end{split}
\end{equation*}
\begin{equation*}
\begin{split}
g(S,S) &= g(Y+\mathbf{\Sigma} W_{j},Y+\mathbf{\Sigma} W_{j}) \\
&=  g_{B}(Y,Y) + \mathbf{\Sigma} (b_{j})^{2} g_{F_{j}}(W_{j},W_{j})
\end{split}
\end{equation*}
\end{theorem}

By taking $Y$ as $h \partial_{t}$ in Theorem \ref{athm}, we
obtain the following result.

\begin{corollary} Let $M=I \times_{b_{1}}  F_{1} \ldots \times_{{b_{m}}}  F_{m}$
be a MGRW space-time with the metric
${\rm g}=-{ \rm d}t^{2} \oplus b_{1}^{2}g_{F_{1}} \oplus \ldots \oplus b_{m}^{2} {\rm g}_{F_{m}}$.
Suppose that
$L=- \partial_{t} +\mathbf{V} $ and $S=h \partial_{t}+\mathbf{W}$ where $\mathbf{V} = \Sigma V_i$
and $\mathbf{W} = \Sigma W_j$ and $h:
I \rightarrow \mathbb{R}$ is a smooth function. Then,

\begin{equation*}
\begin{split}
{ \rm P}(p,\Pi) &=  \mathbf{\Sigma} h b_{k} b_{k}^{\prime \prime}{\rm g}_{F_{k}}(W_{k},V_{k}) \\
&+ \mathbf{\Sigma} h^{2}b_{k} b_{k}^{\prime \prime}{\rm g}_{F_{k}}(V_{k},V_{k}) \\
&- \mathbf{\Sigma} b_{j} b_{j}^{\prime \prime}{\rm g}_{F_{j}}(W_{j},W_{j}) \\
&+ \mathbf{\Sigma}hb_{i} b_{i}^{\prime \prime}{\rm g}_{F_{i}}(W_{i},V_{i}) \\
&- \mathbf{\Sigma}_{j \neq k}   b_{k} (b_{j}^{\prime})^{2}  {\rm g}_{F_{k}}(V_{k},V_{k}) {\rm g}_{F_{j}}(W_{j},W_{j}) \\
&+ \mathbf{\Sigma} b_{i}^{2} {\rm g}_{F_{i}}({\rm R}_{F_{i}}(W_{i},V_{i})V_{i},W_{i}) \\
&- \mathbf{\Sigma}  b_{i}(b_{i}^{\prime})^{4}b_{i}^{\prime \prime}[({\rm g}_{F_{i}}(W_{i},V_{i}))^{2}
-{\rm g}_{F_{i}}(V_{i},V_{i}){\rm g}_{F_{i}}(W_{i},W_{i})]
\end{split}
\end{equation*}
\begin{equation*}
\begin{split}
g(S,S) &= g(h \partial{t}+\mathbf{\Sigma} W_{j},h \partial{t}+\mathbf{\Sigma} W_{j}) \\
&=  -h^{\prime \prime} + \mathbf{\Sigma} (b_{j})^{2} g_{F_{j}}(W_{j},W_{j})
\end{split}
\end{equation*}
\end{corollary}

\subsection{Null Sectional Curvature of a GRW space-time}

Assume that $M=I \times_{b} F$ is a GRW space-time with the metric
$g= g_{I} \oplus b^{2}g_{F}$ where $g_{I}=-dt^{2}$.
Let $\Pi$ be degenerate null plane at $p \in M$ spanned
by a null vector $L$ and space-like vector $S$,i.e.,
$g(L,L)=0$ and $g(S,S)>0$ where $U= \partial_{t}$ is a
reference frame since $g(U,U)= g( \partial_{t}, \partial_{t})=-1$.

\begin{itemize}
\item[(i)] Let  $L= h \partial_{t} + V $, $U= \partial_{t}$

$g(L,U) = g(h \partial_{t} + V , \partial_{t})
= -dt^{2}(h \partial_{t},\partial_{t})
+ b^{2}g_{F}(V,0)\\
= -h$

This implies $g(L,U)=1 $ if and only if $h=-1$. Then we have
$L=-\partial_{t} + V$
\item[(ii)] $g(L,L)=g(- \partial_{t} + V,- \partial_{t} + V)=
-dt^{2}(- \partial_{t},V)+b^{2}g_{F}(V,V)
=-1+b^{2}g_{F}(V,V)
=-1+g_{F}(V,V)$

By imposing $L$ to be null, i.e, $g(L,L)=0$,
we obtain $g_{F}(V,V)=\frac {1}{b^{2}}$.
\item[(iii)]Let  $E =\text{{\rm span}}(\{L,S\})$ be a
degenerate plane section, i.e., $Q(L,S)=0$

$Q(L,S)=Q(- \partial_{t} + V,Y+W)=g(- \partial_{t},- \partial_{t})
g(Y+W,Y+W)
-g(- \partial_{t},Y+W)^{2}
= [- g_{I}(\partial_{t},Y)+ b^{2}g_{F}(V,W)]^{2}
=0$
This implies $g_{F}(V,W)= \frac {1}{b^{2}}g_{I}(\partial_{t},Y)$
\end{itemize}

\begin{corollary}
Let $M=I \times_{b} F$ be a GRW space-time with the metric
$g= -dt^{2} \oplus b^{2}g_{F}$ where $b:I \rightarrow (0, \infty )$
is a smooth function. Then, null sectional curvature of $(M,g)$
is as follows:
\begin{equation}
\begin{split}
g(R(L,S)S,L) &=  -bb^{\prime \prime} g_{F}(W,W) \\
&+ b^{2}g_{F}(R_{F}(W,V)V,W)) \\
&+ bg_{F}(V,V){ \rm H}_{B}^{b}( Y,Y) \\
&+ b^{2} (b^{\prime})^{2}[g_{F}(V,W)^{2}-\frac{1}{b^{2}}g_{F}(W,W)]
\end{split}
\end{equation}

\begin{equation}
\begin{split}
g(S,S) &= g_{I}(Y,Y) + b^{2}g_{F}(W,W)
\end{split}
\end{equation}
\end{corollary}

\begin{remark}
As a special case assume that $Y=0$. Then
\begin{equation}
\begin{split}
{ \rm Q}_{F}(V,W)&=g_{F}(V,V)g_{F}(W,W)-g_{F}(V,W)^{2} \\
&= \frac{1}{b^{2}}g_{F}(W,W)
\end{split}
\end{equation}

In this case we have,

\begin{equation}
\begin{split}
{ \rm K}_{U}(\Pi)&=\frac{1}{b^{2}}K_{F}(V,W)+
\frac{b^{\prime \prime}}{b}-(\frac{b^{\prime}}{b})^{2}
\end{split}
\end{equation}

Using the above equation we establish that ${ \rm K}_{U}(\Pi)=
\frac{1}{b^{2}}{ \rm K}_{F}(V,W)$ if and only if the warping
function $b(t)=ce^{kt}$ where $c$ and $k$ are arbitrary constants.
\end{remark}

\subsection{Null Sectional Curvature of a Generalized Kasner Space-time}

\begin{definition}A generalized Kasner space-time $(M,g)$ is a Lorentzian multiply warped product
of the form  $M=I \times_{\varphi^{p_{1}}}  F_{1} \ldots \times_{\varphi^{p_{m}}}  F_{m}$
with the metric
${\rm g}=-{ \rm d}t^{2} \oplus \varphi^{2p_{1}}g_{F_{1}} \times
\oplus \ldots \oplus \varphi^{2p_{m}} {\rm g}_{F_{m}}$ where $\varphi : I \rightarrow (0, \infty)$ is smooth and $p_{i} \in \mathbb{R}$ for $i=1, \ldots ,m$ and also  $I=(t_{1},t_{2})$ with $ -\infty \leq  t_{1} < t_{2} \leq \infty $.
\end{definition}
\begin{corollary} \label{acorol}
Let $L=- \partial_{t} +\mathbf{V} $ and $S=Y+\mathbf{W}$ where
$\varphi \rightarrow (0, \infty)$ with $- \infty \leq t_{1} < t_{2} \leq \infty $.
and $\mathbf{V} =\mathbf{\Sigma} V_{i}$,  $\mathbf{W} = \mathbf{\Sigma} W_{j}$
Then null sectional of a multiply Kasner space-time of the form above is given by
\begin{equation*}
\begin{split}
{ \rm P}(p,\Pi) &=  \mathbf{\Sigma}\varphi^{p_{k}}{\rm g}_{F_{k}}(W_{k},V_{k})
{\rm H}_{I}^{\varphi^{p_{k}}}(\partial_{t},Y) \\
&+ \mathbf{\Sigma} \varphi^{p_{k}}{\rm g}_{F_{k}}(V_{k},V_{k})
{ \rm H}_{I}^{\varphi^{p_{k}}}(Y,Y) \\
&- \mathbf{\Sigma} \varphi^{p_{j}} {p_{j}}
({p_{j}-1})\varphi^{({p_{j}-2})}{\rm g}_{F_{j}}(W_{j},W_{j}) \\
&+ \mathbf{\Sigma} \varphi^{p_{i}}{\rm g}_{F_{i}}(W_{i},V_{i})
{ \rm H}_{I}^{\varphi^{p_{i}}}(\partial_{t},Y) \\
&- \mathbf{\Sigma}_{j \neq k} \varphi^{p_{k}} {p_{j}}^{2}
\varphi^{2({p_{j}-1})}{\rm g}_{F_{k}}(V_{k},V_{k}) {\rm g}_{F_{j}}(W_{j},W_{j}) \\
&+ \mathbf{\Sigma} \varphi^{2p_{i}} {\rm g}_{F_{i}}({\rm R}_{F_{i}}(W_{i},V_{i})V_{i},W_{i}) \\
&- \mathbf{\Sigma} \varphi^{2p_{i}} p_{i}^{2}\varphi^{2(p_{i}-1)} {p_{i}}({p_{i}}-1)
\varphi^{(p_{i}-2)}    [({\rm g}_{F_{i}}(W_{i},V_{i}))^{2}
-{\rm g}_{F_{i}}(V_{i},V_{i}){\rm g}_{F_{i}}(W_{i},W_{i})]
\end{split}
\end{equation*}

\begin{equation*}
\begin{split}
{\rm g}(S,S) &= {\rm g}(Y+\mathbf{\Sigma} W_{j},Y+\mathbf{\Sigma} W_{j}) \\
&= \mathbf{\Sigma} \varphi^{2p_{j}} {\rm g}_{I}(Y,Y)+ {\rm g}_{F_{j}}(W_{j},W_{j})
\end{split}
\end{equation*}
\end{corollary}

\section{Four-dimensional Space-time Models}

\subsection{Type I: Null Sectional Curvature of a Kasner Space-time
with fiber of dimension $(3)$}

\begin{corollary}Let $M=I \times_{b}F$ and $g= -dt^{2} \oplus  b^{2} g_{F},$
$L= -\partial_{t}+V$ and $S=h \partial_{t}+W$ where $\Pi$ be a degenerate null plane which is spanned by the tangent vectors $L$ and $S$. Then,
\begin{equation} \label{eq9}
\begin{split}
P(\Pi) &= f^{2}bb^{\prime \prime} g_{F}(V,V) \\
&- b b^{\prime \prime}g_{F}(W,W)\\
&+ b b^{\prime \prime}g_{F}(V,W)\\
&+ b^{2} g_{F}(R_{F}(W,V)V,W) \\
&- b^{2}(b^{\prime})^{2}[g_{F}(V,W)^{2}-g_{F}(V,V)g_{F}(V,W)]
\end{split}
\end{equation}
\end{corollary}

\subsection{Type II: Null Sectional Curvature of a Kasner Space-time with fiber of dimension $(1,2)$}

A Kasner space-time with fiber of dimension $(1,2)$ is a Lorentzian multiply warped product $(M,g)$
of the form $M=(0, \infty) \times_{b_{1}} \mathbb{R} \times_{b_{2}} F$
with the metric $g= -dt^{2} \oplus  b_{1}^{2}dx^{2} \oplus b_{2}^{2} g_{F} $

\begin{corollary}Let $\Pi$ be degenerate null plane at $p \in M$ spanned by a null vector
$L$ and space-like vector $S$, i.e., $g(L,L)=0$ and $g(S,S)>0$.

Let $L=- {\partial t} +  f_{1} {\partial x}+ V$ and
$S=f {\partial t} + h_{1} {\partial x}+ W$ where
$$V_{1} = f_{1} {\partial}_{x},  $$
$$W_{1} = h_{1} {\partial}_{x} $$
Then the null
sectional curvature of $\Pi$ with respect to $L$ is is given by
$$K_{U}(\Pi)=\frac{R((L,S)S,L)}{g(S,S)}  $$
where $U= {\partial}_{t}$ is a reference frame since
$g(U,U)= g( {\partial}_{t}, {\partial}_{t})=-1$ and

\begin{equation} \label{}
\begin{split}
g(R(\Sigma V_{i},\Sigma W_{j})\Sigma W_{j},\Sigma V_{i}) &=
b_{1} f_{1}h_{1}f b_{1}^{\prime \prime} + b_{2}f^{2}b_{2}^{\prime \prime}g_{F}(V,V) \\
&- b_{1}  b_{1}^{\prime \prime}h_{1}^{2} -  b_{2} b_{2}^{\prime \prime} g_{F}(W,W)\\
&+ b_{1} f_{1}h_{1}f b_{1}^{\prime \prime} +  b_{2} b_{2}^{\prime \prime} g_{F}(V,W)\\
&+  b_{2}^{2} g_{F}(R_{F}(W,V)V,W)\\
&- b_{2}^{2} (b_{2}^{\prime})^{2}[g_{F}(V,W)^{2} \\
&- g_{F}(V,V)g_{F}(W,W)]
\end{split}
\end{equation}
\begin{equation} \label{eq7a}
\begin{split}
g(S,S) &= g(f {\partial t} +\Sigma W_{j},f{\partial t} +\Sigma W_{j}) \\
&=  -f^{2}+b_{1}^{2}h_{1}^{2}+ b_{2}^{2}g_{F}(W,W)
\end{split}
\end{equation}
\end{corollary}
\begin{remark}
For the case $Y=h \partial_{t}$ and $V_{1}=W_{1}=0$, $P(E)$  becomes as follows:
\begin{equation} \label{eq7b}
\begin{split}
P(\Pi) &= \varphi^{p_{2}} p_{2} (p_{2}-1)\varphi^{p_{2}-2}g_{F_{2}}(W_{2},V_{2}) \\
&+ h^{2} \varphi^{p_{2}}p_{2} (p_{2}-1)\varphi^{p_{2}-2}g_{F_{2}}(V_{2},V_{2}) \\
&- \varphi^{p_{2}} p_{2} (p_{2}-1)\varphi^{p_{2}-2}g_{F_{2}}(W_{2},V_{2}) \\
&+ h \varphi^{p_{2}} p_{2} (p_{2}-1)\varphi^{p_{2}-2}g_{F_{2}}(W_{2},V_{2}) \\
&+ \varphi^{2p_{2}} g_{F_{2}}(R_{F_{2}}(W_{2},V_{2})V_{2},W_{2}) \\
&- p_{2}^{4} \varphi^{6p_{2}-4} [g_{F_{2}}(W_{2},V_{2})^{2}-
g_{F_{2}}(V_{2},V_{2})g_{F_{2}}(W_{2},W_{2})]
\end{split}
\end{equation}
 Similarly we can obtain $P(E)$ for the case $V_{2}=W_{1}=0$ by using the symmetry properties.
\end{remark}
\subsection{Type III: Null Sectional Curvature of a Kasner Space-time with the fiber of dimension $(1,1,1)$}

A Kasner space-time with the fiber of dimension $(1,1,1)$ is a Lorentzian multiply warped product
$(M,g)$ of the form $$M=(0, \infty) \times_{\varphi^{p_{1}}}
\mathbb{R} \times_{\varphi^{p_{2}}} \mathbb{R}
\times_{\varphi^{p_{3}}} \mathbb{R}$$
with the metric

$$g= -dt^{2} \oplus \varphi^{2p_{1}} dx^{2} \oplus
\varphi^{2p_{2}} dy^{2} \oplus \varphi^{2p_{3}} dz^{2}$$

where $p_{1} + p_{2} + p_{3}= p_{1}^{2} + p_{2}^{2} + p_{3}^{2}=1$.

\begin{corollary}Let $\Pi$ be degenerate null plane at $p \in M$ spanned by a null
vector $L$ and spacelike vector $S$,i.e., $g(L,L)=0$ and $g(S,S)>0$.

Let $L=- \partial_{t} + \Sigma V_{i}$ and $S=f \partial_{t} +\Sigma W_{j}$ where
$$V_{1} = f_{1} \partial_{x},V_{2} = f_{2} \partial_{y} $$
$$V_{3} = f_{3} \partial_{z},W_{1} = h_{1} \partial_{x} $$
$$W_{2} = h_{2}\partial_{y},W_{3} = h_{3} \partial_{z} $$

Then the null sectional curvature of $\Pi$ with respect to $L$ is is given by
$$K_{U}(\Pi)=\frac{R((L,S)S,L)}{g(S,S)}  $$

where $U= \partial_{t}$ is a reference frame since $
g(U,U)= g( \partial_{t}, \partial_{t})=-1$.

\begin{equation} \label{}
\begin{split}
g(R(\Sigma V_{i},\Sigma W_{j})\Sigma W_{j},\Sigma V_{i}) &= -\Sigma \varphi^{p_{i}} f_{i} h_{i} \\
&+ \Sigma \varphi^{p_{k}} (f_{k})^{2} f^{2}   p_{k} (p_{k}-1)\varphi^{p_{k}-2} \\
&- \Sigma p_{j} (p_{j}-1)\varphi^{p_{j}-2} (\varphi^{p_{j}}h_{j}^{2}\\
&+ \Sigma \varphi^{p_{i}}f_{i}  h_{i} f  p_{i} (p_{i}-1)\varphi^{p_{i}-2} \\
&- \Sigma_{j \neq k} \varphi^{p_{k}}  f_{k}^{2} p_{j}^{2}\varphi^{(2p_{j}-2)} h_{j}^{2}
\end{split}
\end{equation}
\begin{equation} \label{eq7}
\begin{split}
g(S,S) &= g(f {\partial t} +\Sigma W_{j},f {\partial t} +\Sigma W_{j}) \\
&=  -f^{2} \Sigma\varphi^{2p_{j}} h_{j}^{2}
\end{split}
\end{equation}
\end{corollary}
\section{Null sectional curvature of a SSS-T}

Let $M={f}_I \times F$ be a SSST with the metric ${\rm g}=-f^{2} {\rm g}_{I} \oplus {\rm g}_{F}$
where $-{\rm d}t^{2}$ is the negative definite metric on $I.$

Let $\Pi$ be degenerate null plane at $p \in M$ spanned by a null
vector $L$ and spacelike vector $S$,i.e., ${\rm g}(L,L)=0$ and ${\rm g}(S,S)>0$
with the reference frame  $U=f^{-1} \partial_{t}$ .  The following calculations will be useful in Corollary \ref{sstnull}:
\begin{itemize}
\item[(i)] Let  $L= h \partial_{t} + V $, $U=f^{-1} \partial_{t}$

$g(L,U) = {\rm g}(h \partial_{t} + V ,f^{-1} \partial_{t})
= -f^{2}dt^{2}(h \partial_{t},f^{-1} \partial_{t})
+ {\rm g}_{F}(V,0)\\
= -f^{2}hf^{-1} \\
= -fh$

This implies ${\rm g}(L,U)=1 $ iff $h=-f^{-1}$. Then we have $L=-f^{-1}
\partial_{t} + V$
\item[(ii)] ${\rm g}(L,L)=
{\rm g}(-f^{-1} \partial_{t} + V,-f^{-1} \partial_{t} + V)=
-f^{2}dt^{2}(-f^{-1} \partial_{t},v)+{\rm g}_{F}(V,V)
=-f^{2} \frac {1}{f^{2}}+{\rm g}_{F}(V,V)
=-1+{\rm g}_{F}(V,V)$

If we want to make $g(L,L)=0$, i.e., $L$ is null,
then we must have $g_{F}(V,V)=1$.
\item[(iii)]Let  $E =\text{{\rm span}}(\{L,S\})$
be a degenerate plane section,i.e., $Q(L,S)=0$

$Q(L,S)=Q(-f^{-1} \partial_{t} + V,Y+W)=
{\rm g}(-f^{-1} \partial_{t},-f^{-1} \partial_{t}){\rm g}(Y+W,Y+W)
-{\rm g}(-f^{-1} \partial_{t},Y+W)^{2}
= [-f {\rm g}_{I}(\partial_{t},Y)+ {\rm g}_{F}(V,W)]^{2}
=0$
This implies ${\rm g}_{F}(V,W)=f {\rm g}_{I}(\partial_{t},Y)$
\end{itemize}
\begin{corollary} \label{sstnull}Then the
null sectional curvature of SSST is as follows:
\begin{equation}
\begin{split}
{\rm g}({\rm R}(L,S)S,L) &=  -{\rm g}_{F}(\nabla_{F}f,\nabla_{F}f)
({\rm g}_{I}(Y, \partial_{t})^{2}+{\rm g}_{I}(Y,Y))\\
&- f {\rm g}_{I}(Y,Y){\rm H}_{F}^{f}(V,V) \\
&- {\rm g}_{I}(Y, \partial_{t}){\rm H}_{F}^{f}(V,W) \\
&+ {\rm g}_{I}(Y, \partial_{t}){\rm H}_{F}^{f}(V,W) \\
&+ \frac{1}{f} {\rm H}_{F}^{f}(W,W) \\
&+ {\rm g}_{F}(R_{F}(V,W)V,W)
\end{split}
\end{equation}

If we take $Y= h \partial_{t}$, we obtain that:
\begin{equation}
\begin{split}
{\rm g}({\rm R}(L,S)S,L) &= fh^{2}{\rm H}_{F}^{f}(V,V)\\
&+ \frac{1}{f} {\rm H}_{F}^{f}(W,W)+{\rm g}_{F}(R_{F}(V,W)V,W)
\end{split}
\end{equation}

Since null sectional curvature is independent
from the choice of $S$ (space-like vector),
without loss of generality we can assume that
${\rm g}(S,S)=1$.
Then
\begin{equation}
\begin{split}
{\rm g}({\rm R}(L,S)S,L) &=  -{\rm g}_{F}(\nabla_{F}f,\nabla_{F}f)
({\rm g}_{I}(Y, \partial_{t})^{2}+{\rm g}_{I}(Y,Y))\\
&- f {\rm g}_{I}(Y,Y){\rm H}_{F}^{f}(V,V) \\
&- \frac{1}{f} {\rm H}_{F}^{f}(W,W) \\
&+ {\rm g}_{F}(R_{F}(V,W)V,W)
\end{split}
\end{equation}
\end{corollary}

\begin{remark}
As a special case assume that $Y=0$. Then
\begin{equation}
\begin{split}
Q_{F}(V,W)&={\rm g}_{F}(W,W)
\end{split}
\end{equation}
In this case we have
\begin{equation}
\begin{split}
K_{U}(\Pi)&=K_{F}(V,W)-\frac{H_{F}^{f}(W,W)}{f{\rm g}_{F}(W,W)}
\end{split}
\end{equation}

If we assume that ${\rm H}_{F}^{f}=kf{\rm g}_{F}(W,W)$ for some $k \in \mathbb{R}$, then
$$ K_{U}(\Pi)=K_{F}(V,W)-k$$
\end{remark}
\appendix
\section{Null Sectional Curvature of MGRW space-time}

Here we give main steps of the proof of Theorem \ref{athm}:

\begin{equation*} \label{mgrwlight}
\begin{split}
g(R(L,S)S,L) &= g(R(-\partial_{t} + \mathbf{\Sigma} V_{i},Y+\mathbf{\Sigma} W_{j})Y+
\mathbf{\Sigma} W_{j},- \partial_{t} + \mathbf{\Sigma} V_{i}) \\
&= g(R(-\partial_{t},Y)Y,-\partial_{t}) + g(R(\mathbf{\Sigma} V_{i},Y)Y,-\partial_{t}) \\
&+ g(R(-\partial_{t},\mathbf{\Sigma} W_{j})Y,-\partial_{t}) +
g(R(\mathbf{\Sigma} V_{i},\mathbf{\Sigma} W_{j})Y,-\partial_{t}) \\
&+ g(R(-\partial_{t},Y)\mathbf{\Sigma} W_{j},-\partial_{t})+
g(R(\mathbf{\Sigma} V_{i},Y)\mathbf{\Sigma} W_{j},-\partial_{t}) \\
&+  g(R(-\partial_{t},\mathbf{\Sigma} W_{j})\mathbf{\Sigma} W_{j},-\partial_{t}) +
g(R(\mathbf{\Sigma}V_{i},\mathbf{\Sigma} W_{j})\mathbf{\Sigma} W_{j},-\partial_{t}) \\
&+ g(R(-\partial_{t},Y)Y,\mathbf{\Sigma} V_{i}) +
g(R(\mathbf{\Sigma} V_{i},Y)Y,\mathbf{\Sigma} V_{i}) \\
&+ g(R(-\partial_{t},\mathbf{\Sigma} W_{j})Y,\mathbf{\Sigma} V_{i}) +
g(R(\mathbf{\Sigma} V_{i},\mathbf{\Sigma} W_{j})Y,\mathbf{\Sigma} V_{i}) \\
&+ g(R(-\partial_{t},Y)\mathbf{\Sigma} W_{j},\mathbf{\Sigma} V_{i}) +
g(R(\mathbf{\Sigma} V_{i},Y)\mathbf{\Sigma} W_{j},\mathbf{\Sigma} V_{i}) \\
&+ g(R(-\partial_{t},\mathbf{\Sigma} W_{j})\mathbf{\Sigma} W_{j},\mathbf{\Sigma} V_{i})+
g(R(\mathbf{\Sigma} V_{i},\mathbf{\Sigma} W_{j})\mathbf{\Sigma} W_{j},\mathbf{\Sigma} V_{i})
\end{split}
\end{equation*}

\begin{equation}
\begin{split}
 g(R(-\partial_{t},Y)Y,-\partial_{t})  &= g_{I}(R_{I}(\partial_{t},Y)Y,-\partial_{t}) \\
 &= 0
\end{split}
\end{equation}
\begin{equation}
\begin{split}
 g(R(\mathbf{\Sigma} V_{i},Y)Y,-\partial_{t})  &= \mathbf{\Sigma}
 g(-\frac{1}{b_{i}}{ \rm H}_{I}^{b_{i}}(Y,Y) V_{i},-\partial_{t}) \\
 &=\mathbf{\Sigma}  \frac{1}{b_{i}}{ \rm H}_{B}^{b_{i}}(Y,\partial_{t}) g(V_{i},-\partial_{t}) \\
 &= 0
\end{split}
\end{equation}

\begin{equation} \label{}
\begin{split}
 g(R(-\partial_{t},\mathbf{\Sigma} W_{j})Y,-\partial_{t})  &= \mathbf{\Sigma}
 g(-\frac{1}{b_{j}}{ \rm H}_{I}^{b_{j}}(-\partial_{t},Y) W_{j},-\partial_{t}) \\
 &= \Sigma  \frac{1}{b_{j}}{ \rm H}_{I}^{b_{j}}(\partial_{t},Y) g(W_{j},\partial_{t}) \\
 &= 0
\end{split}
\end{equation}

\begin{equation} \label{}
\begin{split}
 g(R(\mathbf{\Sigma} V_{i},\mathbf{\Sigma} W_{j})Y,-\partial_{t}) &= 0
\end{split}
\end{equation}
\begin{equation} \label{}
\begin{split}
 g(R(-\partial_{t},Y)\mathbf{\Sigma} W_{j},-\partial_{t})  &= 0
\end{split}
\end{equation}
\begin{equation} \label{}
\begin{split}
g(R(\mathbf{\Sigma} V_{i},Y)\mathbf{\Sigma} W_{j},-\partial_{t}) &= \mathbf{\Sigma} g(\frac{-1}{b_{i}}g(V_{i},W_{i})\nabla_{Y}^{I}(\nabla^{I}b_{i}),-\partial_{t})\\
&= \mathbf{\Sigma} \frac{-b_{i}^2}{b_{i}}g_{F_{i}}(V_{i},W_{i})g(\nabla_{Y}^{I}(\nabla^{I}b_{i}),-\partial_{t}) \\
&= \mathbf{\Sigma} b_{i}g_{F_{i}}(V_{i},W_{i}){\rm H}_{I}^{b_{i}}(Y,\partial_{t})
\end{split}
\end{equation}

\begin{equation} \label{}
\begin{split}
g(R(-\partial_{t},\mathbf{\Sigma} W_{j})\mathbf{\Sigma} W_{j},-\partial_{t}) &= \mathbf{\Sigma} g(\frac{-1}{b_{j}}g(V_{i},W_{i})\nabla_{-\partial_{t}}^{I}(\nabla^{I}b_{i}),Y)\\
&= \mathbf{\Sigma} \frac{-1}{b_{i}}g_{F_{j}}(W_{j},W_{j})g(\nabla_{-\partial_{t}}^{I}(\nabla^{I}b_{j}),-\partial_{t}) \\
&= \mathbf{\Sigma} -b_{j}g_{F_{j}}(W_{j},W_{j}){\rm H}_{I}^{b_{i}}(-\partial_{t},-\partial_{t}) \\
&= \mathbf{\Sigma} - b_{j} b_{j}^{\prime \prime} g_{F_{j}}(W_{j},W_{j})
\end{split}
\end{equation}

\begin{equation} \label{}
\begin{split}
g(R(\mathbf{\Sigma} V_{i},\mathbf{\Sigma} W_{j})\mathbf{\Sigma} W_{j},-\partial_{t}) &=
\mathbf{\Sigma}_{i \neq j} g(g(W_{j},W_{j}) \frac{g_{I}(\nabla^{I}b_{j},\nabla^{I}b_{i})}{b_{i}b_{j}}V_{i},\partial_{t}) \\
&+ g(R_{F_{i}}(V_{i},W_{i})W_{i},-\partial_{t})\\
&+ \frac{g_{I}(\nabla^{I}b_{i},\nabla^{I}b_{i})}{b_{i}b_{i}}g(g(V_{i},W_{i})W_{i}-g(W_{i},W_{i})V_{i},-\partial_{t}) \\
&= 0
\end{split}
\end{equation}

\begin{equation} \label{}
\begin{split}
g(R(-\partial_{t},Y)Y,\mathbf{\Sigma} V_{i}) &= g(R_{I}(\partial_{t},Y)Y,\mathbf{\Sigma} V_{i}) \\
&= 0
\end{split}
\end{equation}

\begin{equation} \label{}
\begin{split}
g(R(\mathbf{\Sigma} V_{i},Y)Y,\mathbf{\Sigma} V_{i})  &= \mathbf{\Sigma}
g( -\frac{{\rm H}_{I}^{b_{i}}(Y,-\partial_{t})}{b_{i}}V_{i},V_{i}) \\
 &= \mathbf{\Sigma} b_{i}^{2} \frac{{\rm H}_{B}^{b_{i}}(Y,Y)} {b_{i}} g_{F_{i}}(V_{i},V_{i}) \\
 &= \mathbf{\Sigma} b_{i} {\rm H}_{B}^{b_{i}}(Y,Y)  g_{F_{i}}(V_{i},V_{i})
 \end{split}
\end{equation}

\begin{equation} \label{}
\begin{split}
g(R(-\partial_{t},\mathbf{\Sigma} W_{j})Y,\mathbf{\Sigma} V_{i}) &= \mathbf{\Sigma} g( -\frac{{\rm H}_{I}^{b_{i}}(-\partial_{t},Y)}{b_{i}}W_{j},V_{i}) \\
&= \mathbf{\Sigma}  b_{i}^{2} \frac{{\rm H}_{I}^{b_{i}}(\partial_{t},Y)}{b_{i}} g_{F_{i}}(V_{i},W_{i}) \\
&= \mathbf{\Sigma}  b_{i} {\rm H}_{I}^{b_{i}}(\partial_{t},Y)g_{F_{i}}(V_{i},W_{i})
\end{split}
\end{equation}

\begin{equation} \label{}
\begin{split}
g(R(\mathbf{\Sigma} V_{i},\mathbf{\Sigma} W_{j})Y,\mathbf{\Sigma} V_{i}) &= 0
\end{split}
\end{equation}

\begin{equation} \label{}
\begin{split}
g(R(-\partial_{t},Y)\mathbf{\Sigma} W_{j},\mathbf{\Sigma} V_{i}) &= 0
\end{split}
\end{equation}

\begin{equation} \label{}
\begin{split}
g(R(\mathbf{\Sigma} V_{i},Y)\mathbf{\Sigma} W_{j},\mathbf{\Sigma} V_{i}) &= \mathbf{\Sigma} g(\frac{g(V_{i},W_{i})}{b_{i}}\nabla_{Y}^{I}(\nabla^{I}b_{i}),V_{i}) \\
&= \mathbf{\Sigma} b_{i} g_{F_{i}}(V_{i},W_{i})g(\nabla_{Y}^{I}(\nabla^{I}b_{i}),V_{i})  \\
&= \mathbf{\Sigma} \frac{-g_{F_{i}}(V_{i},W_{i})}{b_{i}} {\rm H}^{I}_{b_{i}}(Y,V_{i}) \\
&= 0
\end{split}
\end{equation}

\begin{equation} \label{}
\begin{split}
g(R(-\partial_{t},\mathbf{\Sigma} W_{j})\mathbf{\Sigma} W_{j},\mathbf{\Sigma} V_{i}) &= \mathbf{\Sigma} g(\frac{g(W_{j},W_{j})}{b_{j}}\nabla_{-\partial_{t}}^{I}(\nabla^{I}b_{i}),V_{i}) \\
&= \mathbf{\Sigma} -b_{j}g_{F_{j}}(W_{j},W_{j}) g(\nabla_{\partial_{t}}^{I}(\nabla^{I}b_{i}),V_{i})) \\
&= 0
\end{split}
\end{equation}

\begin{equation} \label{}
\begin{split}
g(R(\mathbf{\Sigma} V_{i},\mathbf{\Sigma} W_{j})\mathbf{\Sigma} W_{j},\mathbf{\Sigma} V_{i}) &=
\mathbf{\Sigma}_{i \neq j}  g(-g(W_{j},W_{j})  \frac {g_{I}(\nabla^{I}b_{j},\nabla^{I}b_{j})}{b_{j}^2}V_{i},V_{i}) \\
&+ \mathbf{\Sigma} g(R_{F_{i}}(V_{i},W_{i}),W_{i},V_{i}) \\
&+ g(\frac{g_{I}(\nabla^{I}b_{i},\nabla^{I}b_{i})}{b_{i}^2}(g(V_{i},W_{i})W_{i}-g(W_{i},W_{i})V_{i}),V_{i}) \\
&= \mathbf{\Sigma}_{i \neq j} -(b_{j}^{\prime})^2 g_{F_{j}}(W_{j},W_{j})b_{i}^2g_{F_{i}}(V_{i},V_{i}) \\
&+ \mathbf{\Sigma} b_{i}^2g_{F_{i}}(R_{F_{i}}(V_{i},W_{i}),W_{i},V_{i}) \\
&+ \mathbf{\Sigma} b_{i}^2 (b_{i}^{\prime})^2[g_{F_{i}}(V_{i},W_{i})^2-g_{F_{i}}(V_{i},V_{i})g_{F_{i}}(W_{i},W_{i})]
\end{split}
\end{equation}

Then we obtain the following formula for null sectional curvature of Multiply
Generalized Robertson Walker Spacetimes in Theorem \ref{athm}.

\begin{equation*} 
\begin{split}
g(S,S) &= g(Y+\mathbf{\Sigma} W_{j},Y+\mathbf{\Sigma} W_{j}) \\
&= \mathbf{\Sigma} (b_{j})^{2} g_{B}(Y,Y) g_{F_{j}}(W_{j},W_{j})
\end{split}
\end{equation*}





\end{document}